\documentclass{elsarticle}
\usepackage{amssymb,wrapfig,amsmath, amsthm, amscd,ifthen}
\usepackage[cp1251]{inputenc}
\usepackage{color} 
\usepackage{hyperref} 
\usepackage{graphicx,caption,subcaption,mathrsfs,appendix}
\newtheorem{theorem}{Theorem}
\newtheorem{lemma}[theorem]{Lemma}

\newcommand{\1}{\mathbf 1}

\begin{document}
\title{MSO 0-1 law for recursive random trees\tnoteref{t1}}
\tnotetext[t1]{The work of M.E. Zhukovskii is supported by the Ministry of Science and Higher Education of the Russian Federation (Goszadaniye No. 075-00337-20-03), project No. 0714-2020-0005 (the work on Section 2) and by Grant N NSh-2540.2020.1 to support leading scientific schools of Russia (the work on Section 3.1). M.E. Zhukovskii proved results for uniform attachment. The part of the study made by Y.A. Malyshkin was funded by RFBR, project number 19-31-60021. Y.A. Malyshkin proved results for preferential attachment.}
\author[1]{Y.A. Malyshkin}
\ead{yury.malyshkin@mail.ru}
\author[2]{M.E. Zhukovskii}
\ead{zhukmax@gmail.com}

\address[1]{Moscow Institute of Physics and Technology}
\address[2]{Moscow Institute of Physics and Technology}

\begin{abstract}
We prove the monadic second-order 0-1 law for two recursive tree models: uniform attachment tree and preferential attachment tree. We also show that the first order 0-1 law does not hold for non-tree uniform attachment models.
\end{abstract}

\begin{keyword}
uniform attachment tree; preferential attachment tree; 0-1 law; monadic second-order logic
\end{keyword}

\maketitle

\section{Introduction}

Let $n\in\mathbb{N}$. {\it A random graph $\mathcal{G}_n$} is a random element of the set of all undirected graphs without loops and multiple edges on the vertex set $[n]:=\{1,\ldots,n\}$ with a probability distribution $\mu_n$. The case of the uniform distribution $\mu_n$ is widely studied as a particular case of {\it the binomial random graph} denoted by $G(n,p)$~\cite{Bollobas,Janson} where every edge appears independently with probability $p$ (i.e., $\mu_n(G)=p^{|E(G)|}(1-p)^{{n\choose 2}-|E(G)|}$ for every graph $G$ on vertex set $[n]$). Hereinafter, we denote by $V(G)$ and $E(G)$ the set of vertices and the set of edges of $G$ respectively.

Let us recall that {\it a first order (FO) sentence} about graphs expresses a graph property using the following symbols: variables $x,y,x_1,\ldots$, logical connectives $\wedge,\vee,\neg,\Rightarrow,\Leftrightarrow$, predicates $\sim$ (adjacency), $=$ (coincidence), quantifiers $\exists,\forall$ and brackets (see the formal definition in, e.g.,~\cite{Libkin,Survey,Strange}). For example, the property of being complete is expressed by the FO sentence
$$
 \forall x\forall y\quad [\neg(x=y)]\Rightarrow[x\sim y].
$$
A random graph $\mathcal{G}_n$ obeys {\it the FO 0-1 law} if, for every FO sentence $\varphi$, ${\sf P}(\mathcal{G}_n\models\varphi)$ approaches either 0 or 1 as $n\to\infty$. Following traditions of model theory, we write $G\models\varphi$ when $\varphi$ is true on $G$. Study of 0-1 laws for random graph models is closely related to questions about expressive power of formal logics which, in turn, have applications in complexity~\cite{Libkin,Verb}. In 1969 Glebskii, Kogan, Liogon’kii, Talanov~\cite{Glebsk} (and independently Fagin in 1976~\cite{Fagin}) proved that $G(n,\frac{1}{2})$ (i.e., $\mu_n$ is uniform) obeys the FO 0-1 law. In~\cite{Spencer_Ehren}, Spencer proved that, for $p=p(n)$ such that, for every $\alpha>0$, $\min\{p,1-p\}n^{\alpha}\to\infty$ as $n\to\infty$, $G(n,p)$ obeys the FO 0-1 law as well. The sparse case $p=n^{-\alpha}$, $\alpha>0$, was studied in~\cite{Shelah}.

{\it Monadic second order (MSO) logic} is an extension of the FO logic~\cite[Definition 7.2]{Libkin}. Sentences in this logic are built of the same symbols and, additionally, variable unary predicates $X,Y,X_1,\ldots$. For example, the property of being disconnected is expressed by the MSO sentence
$$
 \exists X\quad\biggl[(\exists x\,X(x))\wedge(\exists x\,\neg X(x))\wedge(\forall x\forall y\,[X(x)\wedge \neg X(y)]\Rightarrow[\neg(x\sim y)])\biggr].
$$
In the same way, $\mathcal{G}_n$ obeys {\it the MSO 0-1 law} if, for every MSO sentence $\varphi$, ${\sf P}(\mathcal{G}_n\models\varphi)$ approaches either 0 or 1 as $n\to\infty$. In 1985~\cite{Kaufmann} Kaufmann and Shelah proved that $G(n,\frac{1}{2})$ does not obey the MSO 0-1 law. The same is true for all other constant $p\in(0,1)$ and $p=n^{-\alpha}$, $\alpha\in(0,1]\cup\{1+1/\ell,\,\ell\in\mathbb{N}\}$ (see~\cite{Zhuk_Ostr_APAL,Tyszkiewicz,Zhuk_JML}).

Further, many other random graph models were studied in the context of logical limit laws. Let us list some of them. In~\cite{McColm}, it was proven that the labeled uniform random tree ($\mu_n(T)=n^{2-n}$ for every tree $T$ on vertex set $[n]$) obeys the MSO 0-1 law (earlier, in~\cite{Woods}, a related result was obtained using generating functions, and the same proof works for this result as well). The FO behavior of random regular graphs was studied in~\cite{Haber}. In~\cite{geo}, logical laws were proven for random geometric graphs. In~\cite{Muller}, the FO and the MSO 0-1 laws were studied for minor-closed classes of graphs. In~\cite{Zhuk_Svesh}, FO 0-1 laws were proven for the classical uniform random graph model $G(n,m)$ ($m$ edges are chosen uniformly at random). Finally, some results related to the FO behavior of the preferential attachment random graph model were obtained in~\cite{Kleinberg}.

In this paper, we study the logical behavior of two well-known recursive random graph models: the uniform model and the preferential attachment model~\cite{recursive}. Let $m\in\mathbb{N}$. The uniform attachment random graph $G^{\mathrm{U}}(n,m)$ is defined recursively: $G^{\mathrm{U}}(m+1,m+1)$ is the complete graph on $[m+1]$; for every $n\geq m+1$, $G^{\mathrm{U}}(n+1,m)$ is obtained from $G^{\mathrm{U}}(n,m)$ by adding the vertex $n+1$ with $m$ edges going from $n+1$ to vertices from $[n]$ chosen uniformly at random:
\begin{align*}
 {\sf P}\biggl(n+1\sim x_1,\ldots,n+1\sim x_m\text{ in }G^{\mathrm{U}}(n+1,m)\biggr)={n\choose m}^{-1},\\ 1\leq x_1<\ldots<x_m\leq n.
\end{align*}
In Section~\ref{no_law}, we show that, for every $m\geq 2$, $G^{\mathrm{U}}(n,m)$ does not obey the FO 0-1 law. For $m=1$, we prove the following positive result.
\begin{theorem}
$G^{\mathrm{U}}(n,1)$ obeys the MSO 0-1 law. 
\end{theorem}

In the preferential attachment random graph $G^{\mathrm{P}}(n,m)$, we also start from the complete graph $G^{\mathrm{P}}(m+1,m+1)$. $G^{\mathrm{P}}(n+1,m)$ is also obtained from $G^{\mathrm{P}}(n,m)$ by adding the vertex $n+1$ with $m$ edges going from $n+1$ to vertices from $[n]$. The only difference is that these edges $e_1,\ldots,e_m$ are drawn independently, each one has distribution ${\sf P}(e_i=\{n+1,v\})=\frac{\mathrm{deg}_{G^{\mathrm{P}}(n,m)}v}{2mn}$, $v\in[n]$. Notice that this graph may have multiple edges in contrast to all the previous. In the context of 0-1 laws in the considered logic, it is convenient to remove all repetitions and consider a simple graph instead. 
Notice that this modification does not change the model when $m=1$. In~\cite{Kleinberg}, it was proven that  $G^{\mathrm{P}}(n,m)$ does not obey the FO 0-1 law for every $m\geq 3$. 
In this paper, we prove that the MSO 0-1 law holds for $m=1$.
\begin{theorem}
$G^{\mathrm{P}}(n,1)$ obeys the MSO 0-1 law. 
\end{theorem}

The question about validity of both the FO and the MSO 0-1 law for $G^{\mathrm{P}}(n,2)$ remains open.

\section{FO 0-1 law fails for the uniform model when $m\geq 2$} 
\label{no_law}

Let us first assume that $m=2$. Let $X_n$ be the number of {\it diamond graphs} (graph with 4 vertices and 5 edges) in $G^{\mathrm{U}}(n,2)$. 

If $G^{\mathrm{U}}(n,2)$ contains a diamond graph on vertices $u_1<u_2<u_3<u_4$, then  $u_1$, $u_2$ and $u_3$ are adjacent, and $u_4$ is adjacent to exactly two of $u_1,u_2,u_3$. So, there are three ways of distributing edges among $u_1,u_2,u_3,u_4$ to get a diamond graph. Moreover, there are exactly $2u_3-5$ edges in $G^{\mathrm{U}}(u_3-1,2)$ (the process starts from 3 vertices and 3 edges and then, at each step, $1$ vertex and $2$ edges are introduced). Hence, the probability that the neighbors of $u_3$ in $G^{\mathrm{U}}(u_3,2)$ are adjacent is exactly $\frac{2u_3-5}{{{u_3-1}\choose 2}}$.

Let us fix vertices $u_3<u_4$. The probability of the existence of $u_1,u_2$ such that $u_1<u_2<u_3$ and vertices $u_1,u_2,u_3,u_4$ induces a diamond graph in $G^{\mathrm{U}}(u_4,2)$ equals the product of $\frac{2u_3-5}{{{u_3-1}\choose 2}}$ and the probability that all neighbors of $u_4$ in $G^U(u_4,2)$ are among $u_3$ and its neighbors in $G^{U}(u_3,2)$. So, this probability equals $\frac{3(2u_3-5)}{{{u_3-1}\choose 2}{{u_4-1}\choose 2}}$.


Hence,
$$
{\sf E}X_n=\!\sum_{3\leq u_3<u_4\leq n}\frac{3(2u_3-5)}{{{u_3-1}\choose 2}{{u_4-1}\choose 2}}=\!\sum_{3\leq u_3<u_4\leq n}\frac{12(2u_3-5)}{(u_3\!-\!1)(u_3\!-\!2)(u_4\!-\!1)(u_4\!-\!2)}\to\beta
$$
where $\beta>0$ is finite. The latter is immediate due to the integral test of convergence. It also can be verified in the following elementary way:
\begin{align*} 
{\sf E}X_n<\quad24 &\sum_{i=1}^{n-1}\frac{1}{i}\sum_{j=i+1}^{n-1} \frac{1}{j(j+1)}\\
=&\sum_{i=1}^n \frac{24}{i}\left(\frac{1}{i+1}-\frac{1}{n}\right)<
\sum_{i=1}^n \frac{24}{i(i+1)}=24\left(1-\frac{1}{n+1}\right)<24.
\end{align*}

For every $k>3$, denote by $g(k)$ the maximum value of $X_k$, i.e. ${\sf P}(X_k=g(k))>0$ while ${\sf P}(X_k>g(k))=0$. At time $k$, a new diamond graph appears if the new vertex is adjacent to vertices $u,v$ such that $u$ is adjacent to $v$ and there exists a triangle containing $u,v$ in $G^{\mathrm{U}}(k-1,2)$. The number of new diamond graphs equals to the number of such triangles. Therefore, at time $k$ appears at most $k-3$ new triangles. Thus, $g(k)\leq {{k-2}\choose 2}$. Since, with positive probability, vertices $3,\ldots,k$ of $G^{\mathrm{U}}(k,2)$ are adjacent to both 1 and 2, and $X_k={{k-2}\choose 2}$ on such graph, the upper bound is achievable, i.e. $g(k)={{k-2}\choose 2}$. 

Let us consider a FO sentence $\varphi_k$ describing the property of having at least $g(k)$ diamond graphs and prove that, for $k$ large enough, ${\sf P}(G^{\mathrm{U}}(n,2)\models\varphi_k)$ does not approach neither 0 nor 1.

Let us first prove that, for any $k$, ${\sf P}(G^{\mathrm{U}}(n,2)\models\varphi_k)$ is bounded away from 0. Since $X_n$ is non-decreasing function of $n$ a.s., for $n>k$ we get
\begin{equation}
{\sf P}(X_n\geq g(k))\geq{\sf P}(X_k=g(k))>0.
\label{non_conv_below}
\end{equation}

Now, let us prove that, for $k$ large enough, ${\sf P}(G^{\mathrm{U}}(n,2)\models\varphi_k)$ is bounded away from 1.
Fix $\varepsilon>0$ and choose $k$ in a way such that $\frac{\beta}{g(k)}<1-\varepsilon$. Then, by Markov's inequality, for $n$ large enough (so that ${\sf E}X_n<\beta\frac{1-\varepsilon/2}{1-\varepsilon}$),
\begin{equation}
{\sf P}(X_n\geq g(k))\leq\frac{{\sf E}X_n}{g(k)}<1-\frac{\varepsilon}{2}.
\label{non_conv_above}
\end{equation}
Having this, we conclude that $G^{\mathrm{U}}(n,2)$ does not obey the FO 0-1 law.\\

Finally, let $m\geq 3$. Let $X_n$ be the number of $K_{m+1}$ (complete graphs on $m+1$ vertices) in $G^{\mathrm{U}}(n,m)$. 

Let $1\leq u_1<\ldots<u_m<u_{m+1}\leq n$. Then, for $j\in\{1,\ldots,m\}$, the probability that $u_{j+1}$ is adjacent to $u_1,\ldots,u_j$ in $G^{\mathrm{U}}(n,m)$ equals $\frac{{{u_{j+1}-1-j}\choose {m-j}}}{{{u_{j+1}-1}\choose m}}$ (here, we set ${i\choose j}:=1$ when $0\leq i<j$). Therefore,
$$
{\sf P}\biggl(G^{\mathrm{U}}(n,m)|_{\{u_1,\ldots,u_{m+1}\}}\text{ is a clique}\biggr)=\frac{{{u_2-2}\choose {m-1}}}{{{u_2-1}\choose m}}\frac{{{u_3-3}\choose {m-2}}}{{{u_3-1}\choose m}}\ldots
\frac{{{u_m-m}\choose 1}}{{{u_m-1}\choose m}}\frac{1}{{{u_{m+1}-1}\choose m}}
$$
(hereinafter, we denote by $H_U$ a subgraph of $H$ induced by $U$).
We get
\begin{equation*} 
\begin{split}
{\sf E}X_n & = \sum_{1\leq u_1<\ldots<u_{m+1}\leq n}\prod_{j=1}^m\frac{{{u_{j+1}-1-j}\choose {m-j}}}{{{u_{j+1}-1}\choose m}} \\
 & = \sum_{k=1}^m{m\choose k}\sum_{m+1\leq u_{k+1}<\ldots<u_{m+1}\leq n}\prod_{j=k}^m\frac{{{u_{j+1}-1-j}\choose {m-j}}}{{{u_{j+1}-1}\choose m}} \\
 & + \sum_{m+1\leq u_1<\ldots<u_{m+1}\leq n}\prod_{j=1}^m\frac{{{u_{j+1}-1-j}\choose {m-j}}}{{{u_{j+1}-1}\choose m}},
\end{split}
\end{equation*}
where $k$ denotes the index of the maximum label among $u_1,\ldots,u_{m+1}$ that less then $m+1$ (i.e., $u_k\leq m<u_{k+1}$). Notice that in the last equality, we divide the summation into two parts --- in the first part, we sum up over those $u_1,\ldots,u_m$ that contain values at most $m$, i.e. such that $u_1\leq m$.


Finally, we bound ${\sf E}X_n$ from above by
$$
 \sum_{k=1}^m\!\!{m\choose k}\!\!\sum_{m+1\leq u_{k\!+\!1}<\ldots<u_{m\!+\!1}\leq n}\prod_{j=k}^m\frac{m^j}{(u_{j\!+\!1}\!-\!j)^j}\!+\!\!\sum_{m+1\leq u_1<\ldots<u_{m\!+\!1}\leq n}\prod_{j=1}^m\frac{m^j}{(u_{j\!+\!1}\!-\!j)^j}.
$$
Since $m>2$, this upper bound converges to a finite limit, which could be verified using the integral test for convergence (the second term converges since the integral 
$\int_{m+1}^{\infty}\int_{x_1}^{\infty}\ldots\int_{x_m}^{\infty}\left(\prod_{j=1}^m(x_{j+1}-j)^j\right)^{-1}dx_1...dx_{m+1}$
 converges if $\sum_{j=1}^m j>m+1$, which is true for $m>2$, and the same applies to the first summation). Since ${\sf E}X_n$ increases in $n$, we get that, for some $\beta$, ${\sf E}X_n\to\beta$ as $n\to\infty.$

The rest of the proof is the same as in the case $m=2$. For $k>m+1$, $g(k)=k-m$ is the maximum value of $X_k$. Here, we consider a FO sentence $\varphi_k$ describing the property of having at least $g(k)$ $(m+1)$-cliques. Then, we
choose $k$ in a way such that $\frac{\beta}{g(k)}<1-\varepsilon$. In the same way, relations (\ref{non_conv_below})~and~(\ref{non_conv_above}) hold. Therefore, $G^{\mathrm{U}}(n,m)$ does not obey the FO 0-1 law.\\

\section{Proofs}

For a tree $G$ and a vertex $R\in V(G)$, we denote by $G_R$ the tree $G$ rooted in $R$, i.e. the tree with the defined parent-children relation. Rooted trees $G_u$ and $H_v$ are {\it isomorphic} (denoted by $G_u\cong H_v$) if there exists a bijection $f:V(G_u)\to V(H_v)$ that preserves the child--parent relation: $a$ is a child of $b$ in $G_u$ if and only if $f(a)$ is a child of $f(b)$ in $H_v$.

Given a tree $\mathcal{T}$ and a rooted tree $G_R$, we say that {\it $\mathcal{T}$ has a pendant $G_R$}, if there is an edge $\{u,v\}$ in $\mathcal{T}$ such that, after its deletion, the component $F$ of $\mathcal{T}$ containing $v$ is such that $F_v\cong G_R$.\\

We will use the following claim proved in~\cite{McColm} (hereinafter, given a graph property $P$, we say that $\mathcal{G}_n$ has $P$ {\it with high probability}, if $\mu_n(P)\to 1$ as $n\to\infty$).

\begin{lemma}{\cite[Theorem 2.1]{McColm}}
Let $\mathcal{G}_n$ be a random tree (i.e. $\mu_n$ is positive only on trees). For every rooted tree $G_R$, suppose that with high probability $\mathcal{G}_n$ has a pendant $G_R$. Then $\mathcal{G}_n$ obeys the MSO 0-1 law.
\label{the_tool}
\end{lemma}

\subsection{MSO 0-1 law for the uniform recursive tree}

By Lemma~\ref{the_tool}, it is sufficient to prove that, for every rooted $G_R$, $G^{\mathrm{U}}(n,1)$  contains a pendant $G_R$ with high probability.

Consider an arbitrary rooted tree $G_R$. Let $v$ be the number of vertices of $G_R$. Let $R=j_1<\ldots<j_v$ be a labelling of vertices of $G_R$ such that, for every $s\in\{2,\ldots,v\}$, $j_s$ is adjacent to $j_{i}$ for some $i<s$. 

Let $n_0,r\in\mathbb{N}$ and let $n\geq n_0+r+v$. Let $n_0+r<i_1<\ldots<i_v\leq n$. Let $B_{i_1,...,i_v}(n_0,r,n)$ denote the event that, in $G^{\mathrm{U}}(n,1)$, there is an edge between $i_1$ and $n_0$ and its deletion divides the tree into two connected components such that one of them (denote it by $H$) consists of $i_1,\ldots,i_v$ and the bijection $j_s\to i_s$, $s\in\{1,\ldots,v\}$, is an isomorphism of $G_R$ and $H_{i_1}$ (tree $H$ with the root in $i_1$). Let 
$$
\tilde B_{i_1,\ldots,i_v}:=\tilde B_{i_1,\ldots,i_v}(n_0,n,r)=\bigsqcup_{\ell=0}^{r-1}B_{i_1,\ldots,i_v}(n_0+\ell,r-\ell,n),
$$
$$
X:=X(n_0,n,r)=\sum_{n_0+r<i_1<\ldots<i_v\leq n}I_{\tilde B_{i_1,...,i_v}(n_0,n,r)}.
$$
Notice that the event $\{X>0\}$ implies the existence of a pendant $G_R$ in $G^{\mathrm{U}}(n,1)$. So, it is sufficient to prove that for every $\varepsilon>0$ there exists $r\in\mathbb{N}$ such that ${\sf P}(X>0)>1-\varepsilon$ for all large enough $n$.

We get that for $\ell<r$ (we put $i_{v+1}:=n+1$)
\begin{equation}
\begin{aligned}
 {\sf P}(B_{i_1,...,i_v}(n_0+\ell,r-\ell,n))& =\prod_{s=1}^{v}\left(\frac{1}{i_s-1}\prod_{t=i_s+1}^{i_{s+1}-1}\left(1-\frac{s}{t-1}\right)\right)\\
& =\frac{\prod_{s=1}^{v}\left(\prod_{t=i_s+1}^{i_{s+1}-1}\left(t-1-s\right)\right)} {\prod_{s=1}^{v}\left((i_s-1)\prod_{t=i_s+1}^{i_{s+1}-1}\left(t-1\right)\right)}\\
=\frac{(n-1-v)!/(i_1-2)!}{(n-1)!/(i_1-2)!}& =\frac{(n-1-v)!}{(n-1)!}=\frac{1}{(n-v)\ldots(n-1)}.
\end{aligned}
\label{not_depend_on_i}
\end{equation}
(the third equality follows from that fact that, for every $s$,  $(i_{s+1}-1)-1-s$ and $(i_{s+1}+1)-1-(s+1)$ are consecutive numbers).
Since, for different $\ell$, the events $B_{i_1,...,i_v}(n_0+\ell,r-\ell,n)$ are disjoint, we get
$$ 
{\sf P}(\tilde B_{i_1,...,i_v})=\frac{r}{(n-v)\ldots(n-1)}.
$$
Therefore, for any $n_0$
$$
 {\sf E}X={n-n_0-r\choose v}\frac{r}{(n-v)\ldots(n-1)}\to \frac{r}{v!},\quad n\to\infty.
$$

For distinct sets $(i_1,\ldots,i_v)$ and $(\tilde i_1,\ldots,\tilde i_v)$, the events $\tilde B_{i_1,...,i_v}$ and $\tilde B_{\tilde i_1,...,\tilde i_v}$ are disjoint if $\{i_1,\ldots,i_v\}\cap\{\tilde i_1,\ldots,\tilde i_v\}\neq\varnothing$ (in particular, it immediately implies that $\mathrm{cov}\left(I_{\tilde B_{i_1,...,i_v}},I_{\tilde B_{\tilde i_1,...,\tilde i_v}}\right)<0$). Otherwise (if $i$ and $\tilde{i}$ are disjoint), let  $(\sigma_1,\ldots,\sigma_{2v})$ be the permutation of $(i_1,\ldots,i_v,\tilde i_1,\ldots,\tilde i_v)$ such that $\sigma_1<\ldots<\sigma_{2v}$. Then, letting $\sigma_{2v+1}:=n+1$, in the same way as in~(\ref{not_depend_on_i}) we get
\begin{align*}
{\sf P}(\tilde B_{i_1,...,i_v}(n_0,n,r)\cap \tilde B_{\tilde i_1,...,\tilde i_v}(n_0,n,r))\\
=\sum_{\ell=0}^{r-1}\sum_{k=0}^{r-1}{\sf P}(B_{i_1,...,i_v}(n_0+\ell,r-\ell,n)\cap B_{\tilde i_1,...,\tilde i_v}(n_0+k,r-k,n))\\
=r^2
\prod_{s=1}^{2v}\left(\frac{1}{\sigma_s-1}\prod_{t=\sigma_s+1}^{\sigma_{s+1}-1}\left(1-\frac{s}{t-1}\right)\right)
=\frac{r^2}{(n-2v)\ldots(n-1)}.
\end{align*}
Therefore, we get
\begin{align*} 
\mathrm{Var}X&=\mathrm{Var}\left(\sum_{n_0+r<i_1<\ldots<i_v\leq n}I_{\tilde B_{i_1,...,i_v}}\right)\\
&=\sum_{n_0+r<i_1<\ldots<i_v\leq n}\mathrm{Var}\left[I_{\tilde B_{i_1,...,i_v}}\right]\\
&+ \sum_{{\scriptsize \begin{array}{c}n_0+r<i_1<\ldots<i_v\leq n,\\n_0+r<\tilde i_1<\ldots<\tilde i_v\leq n, \\(i_1,\ldots,i_v)\neq(\tilde i_1,\ldots,\tilde i_v)\end{array}}}
\mathrm{cov}\left(I_{\tilde B_{i_1,...,i_v}},I_{\tilde B_{\tilde i_1,...,\tilde i_v}}\right)\\
&<\sum_{n_0+r<i_1<\ldots<i_v\leq n}\left({\sf P}_{\tilde B_{i_1,...,i_v}}-\left[{\sf P}_{\tilde B_{i_1,...,i_v}}\right]^2\right)\\
&+\!\!\!\sum_{{\scriptsize \begin{array}{c}n_0+r<i_1<\ldots<i_v\leq n,\\n_0+r<\tilde i_1<\ldots<\tilde i_v\leq n, \\(i_1,\ldots,i_v)\!\cap\!(\tilde i_1,\ldots,\tilde i_v)=\varnothing\end{array}}}\!\!\! \left[{\sf P}\!\left(\tilde B_{i_1,...,i_v}\!\cap\!\tilde B_{\tilde i_1,...,\tilde i_v}\right)\!-\!{\sf P}\!\left(\tilde B_{i_1,...,i_v}\right)\!{\sf P}\!\left(\tilde B_{\tilde i_1,...,\tilde i_v}\right)\right]\\
&\!=\!{n\!-\!n_0\!-\!r\choose v}\!\left(\frac{r}{(n-v)\ldots(n-1)}-\left(\frac{r}{(n-v)\ldots(n-1)}\right)^2\right)\\
 &\!+\!{n\!-\!n_0\!-\!r\choose v}\!{n\!-\!n_0\!-\!r\!-\!v\choose v}\!\left(\!\frac{r^2}{(n\!-\!2v)\!\ldots\!(n\!-\!1)}\!-\!\left(\frac{r}{(n\!-\!v)\!\ldots\!(n\!-\!1)}\right)^2\right)\\
 &\!=\!{n\!-\!n_0\!-\!r\choose v}\frac{r}{(n-v)\ldots(n-1)}+O\left(\frac{1}{n}\right)\to\frac{r}{v!},\quad\quad\quad\quad n\to\infty.
\end{align*}
It remains to apply Chebyshev's inequality:
$$
 {\sf P}(X=0)\leq\frac{\mathrm{Var}X}{({\sf E}X)^2}\to\frac{v!}{r},\quad n\to\infty.
$$

\subsection{MSO 0-1 law for the preferential attachment random tree}

As above, here, we prove that, for every rooted $G_R$, $G^{\mathrm{P}}(n,1)$  contains a pendant $G_R$ with high probability.

In the same way, we consider a labelling $R=j_1<\ldots<j_v$ of vertices of $G_R$ such that, for every $s\in\{2,\ldots,v\}$, $j_s$ is adjacent to $j_{i}$ for some $i<s$.  

Let $i_1<\ldots<i_v\leq n$. Let $B_{i_1,...,i_v}(n)$ denote the event that, in $G^{\mathrm{P}}(n,1)$, there exists a vertex $n_0<i_1$ adjacent to $i_1$ such that deletion of the edge $\{n_0,i_1\}$  divides the tree into two connected components $H$ and $G^{\mathrm{P}}(n,1)\setminus H$ such that $H$ is induced by $i_1,\ldots,i_v$ and the bijection $j_s\to i_s$, $s\in\{1,\ldots,v\}$, is an isomorphism of $G_R$ and $H_{i_1}$. Let 
$$
X=X(n)=\sum_{2\leq i_1<\ldots<i_v\leq n}I_{B_{i_1,...,i_v}(n)}.
$$
As above, the event $\{X>0\}$ implies the existence of a pendant $G_R$ in $G^{\mathrm{P}}(n,1)$. 

Notice that, for every $\nu\in\{1,\ldots,v\}$, and $i_\nu<s<i_{\nu+1}$ (hereinafter, $i_{v+1}=n+1$), under the condition that, in $G^{\mathrm{P}}(s-1,1)$, deletion of the edge from $i_1$ to an older vertex 
separates vertices $i_1,\ldots,i_\nu$ from the rest of the tree, the probability that $s$ is not adjacent to any of $i_1,\ldots,i_\nu$ in $G^{\mathrm{P}}(n,1)$ equals $1-\frac{2\nu-1}{2(s-1)}$. 

For $\ell\in\{2,\ldots,v\}$, let $x_\ell<j_{\ell}$ be the neighbor of $j_{\ell}$ in the induced subgraph $G_R|_{\{j_1,\ldots,j_{\ell}\}}$. Denote $d_\ell\!=\!\mathrm{deg}_{G_R|_{\{j_1,\ldots,j_{\ell-1}\}}}x_{\ell}$ if $x_{\ell}\!\neq\! j_1$ and   $d_\ell\!=\!\mathrm{deg}_{G_R|_{\{j_1,\ldots,j_{\ell-1}\}}}x_{\ell}+1$ if $x_{\ell}=j_1$. Set $D:=\prod_{\ell=2}^{v}d_{\ell}$. As usual, we denote $(2k+1)!!:=\prod_{i=1}^k(2i+1)=\frac{(2k+1)!}{k!2^k}$.  Then, using Stirling's formula, 
 we get
$$ 
{\sf P}(B_{i_1,...,i_v}(n)) =\left(\prod_{t=i_1+1}^{i_2-1}\left(1-\frac{1}{2(t-1)}\right)\right) \prod_{\ell=2}^{v}\left(\frac{d_{\ell}}{2(i_{\ell}-1)} \prod_{t=i_{\ell}+1}^{i_{\ell+1}-1}\left(1-\frac{2\ell-1}{2(t-1)}\right)\right)
$$
$$
=D\frac{\left(\prod_{t=i_1+1}^{i_2-1}(2t-3)\right)\left(\prod_{\ell=2}^{v}\left(\prod_{t=i_{\ell}+1}^{i_{\ell+1}-1}(2t-2\ell-1)\right)\right)} {\left(\prod_{t=i_1+1}^{i_2-1}(2(t-1))\right)\left(\prod_{\ell=2}^{v}2(i_{\ell}-1)\left(\prod_{t=i_{\ell}+1}^{i_{\ell+1}-1}(2(t-1))\right)\right)}
$$
$$
=D\frac{\prod_{t=i_1}^{n-v}(2t-1)}{2^{n-i_1}\prod_{t=i_1+1}^{n}(t-1)}
=D\frac{[2(n-v)-1]!!/[2(i_1-1)-1]!!}{2^{n-i_1}[n-1]!/[i_1-1]!} 
$$
$$ =D\frac{\frac{(2(n-v))!2^{i_1-1}(i_1-1)!}{2^{n-v}(n-v)!(2(i_1-1))!}}{2^{n-i_1}[n-1]!/[i_1-1]!}
=D\frac{2^{2i_1-2}(i_1-1)!(i_1-1)!}{(2(i_1-1))!}\frac{(2(n-v))!}{2^{2n-v-1}(n-v)!(n-1)!}
$$
$$
=D\frac{2^{2i_1-2}\sqrt{2\pi (i_1-1)}(i_1-1)^{i_1-1}\sqrt{2\pi (i_1-1)}(i_1-1)^{i_1-1}e^{2(i_1-1)-(i_1-1)-(i_1-1)}}{\sqrt{2\pi (2(i_1-1))}(2(i_1-1))^{2(i_1-1)}}
$$
$$
\times\frac{\sqrt{2\pi (2n-2v)}(2(n-v))^{2(n-v)}e^{n-v+n-1-2(n-v)}}{2^{2n-v-1}\sqrt{2\pi (n-v)}(n-v)^{n-v}\sqrt{2\pi (n-1)}(n-1)^{n-1}}\left(1+O\left(\frac{1}{i_1}\right)\right)
$$
$$
=D\sqrt{(i_1-1)} \frac{e^{v-1}\left(1-\frac{v-1}{n-1}\right)^{n-v}}{2^{v-1}(n-1)^{v-1/2}} \left(1+O\left(\frac{1}{i_1}\right)\right)
=\frac{D}{2^{v-1}}\sqrt{\frac{i_1}{n^{2v-1}}}\left(1+O\left(\frac{1}{i_1}\right)\right).
$$
Therefore,
\begin{align*} 
 {\sf E}X&=\sum_{1\leq i_1<\ldots<i_v\leq n}{\sf P}\left(B_{i_1,...,i_v}(n)\right)\\
&=\frac{D}{2^{v-1}}\sum_{i_1=1}^{n-v+1}\sqrt{\frac{i_1}{n^{2v-1}}}{n-i_1\choose v-1}\left(1+O\left(\frac{1}{i_1}\right)\right)\\
&=\frac{D}{(v-1)!2^{v-1}}\sum_{i_1=1}^{n-v+1}\sqrt{\frac{i_1}{n}}\left(1-\frac{i_1}{n}\right)^{v-1}\left(1+O\left(\frac{1}{i_1}\right)+O\left(\frac{1}{n-i_1}\right)\right)\\
 &\sim\frac{2Dn}{(2v+1)!!},\quad n\to\infty.
\end{align*}

The latter relation is derived from the following approximations: 
\begin{align*}
 \sum_{i_1=1}^{n-v+1}\!\!\sqrt{\frac{i_1}{n}}\!\left(\!1\!-\!\frac{i_1}{n}\!\right)^{v-1}\!\!&=n\int_0^1\sqrt{x}(1-x)^{v-1}dx+O(1)\\
 &=n\mathrm{B}\left(\frac{3}{2},v\right)+O(1)= n\frac{\Gamma(\frac{3}{2})\Gamma(v)}{\Gamma(\frac{3}{2}+v)}+O(1)\\
 &=n\frac{\frac{1}{2}\Gamma(\frac{1}{2})(v\!-\!1)!}{(v\!+\!\frac{1}{2})(v\!-\!\frac{1}{2})\ldots\frac{1}{2}\Gamma(\frac{1}{2})}\!+\!O(1)\!=\!n\frac{2^v(v\!-\!1)!}{(2v\!+\!1)!!}\!+\!O(1),
\end{align*}
\begin{align*}
 \sum_{i_1=1}^{n-v+1}\sqrt{\frac{1}{ni_1}}\left(1-\frac{i_1}{n}\right)^{v-1}=& \int_0^1\frac{(1-x)^{v-1}}{\sqrt{x}}dx+O\left(\frac{1}{n}\right)\\
 &=\mathrm{B}\left(\frac{1}{2},v\right)+O\left(\frac{1}{n}\right)=O(1),
 \end{align*}
  \begin{align*}
 \sum_{i_1=1}^{n-v+1}\sqrt{\frac{i_1}{n(n-i_1)}}\left(1-\frac{i_1}{n}\right)^{v-1}&\\
 =\frac{1}{n}\sum_{i_1=1}^{n-v+1}\sqrt{\frac{i_1}{n}}\left(1-\frac{i_1}{n}\right)^{v-2}=&
 \mathrm{B}\left(\frac{3}{2},v-1\right)+O\left(\frac{1}{n}\right)=O(1).\\
\end{align*}

For distinct sets $(i_1,\ldots,i_v)$, $(\tilde i_1,\ldots,\tilde i_v)$, the events $B_{i_1,...,i_v}(n)$ and $B_{\tilde i_1,...,\tilde i_v}(n)$ are disjoint if $\{i_1,\ldots,i_v\}\cap\{\tilde i_1,\ldots,\tilde i_v\}\neq\varnothing$ (in particular, it immediately implies that $\mathrm{cov}\left(I_{B_{i_1,...,i_v}}(n),I_{B_{\tilde i_1,...,\tilde i_v}}(n)\right)<0$). Otherwise, assume that $i_1<\tilde i_1$ and let $\mu\in\{1,\ldots,v\}$ be such that $i_{\mu}<\tilde i_1<i_{\mu+1}$. Let $(\sigma_1,\ldots,\sigma_{2v})$ be the permutation of $(i_{1},\ldots,i_v,\tilde i_1,\ldots,\tilde i_v)$ such that $\sigma_1<\ldots<\sigma_{2v}$. Let $\sigma_{2v+1}=n+1$. Then, for $s\in\{1,\ldots,2v\}$, under the condition that, in $G^{\mathrm{P}}(\sigma_s,1)$, deletion of the edge from $i_1$ to its older neighbor separates vertices $\sigma_1,\ldots,\sigma_s$ from the rest of the tree, the probability that a vertex $t\in(\sigma_s,\sigma_{s+1})$ is not adjacent to any of $\sigma_1,\ldots,\sigma_s$ in $G^{\mathrm{P}}(n,1)$ equals $1-\frac{2\sigma_s-1}{2(t-1)}$ if $s\leq\mu$ and equals $1-\frac{2\sigma_s-2}{2(t-1)}$ if $s>\mu$. Hence, computing the joint probability in the same way as ${\sf P}(B_{i_1,...,i_v}(n))$, we get
$$
{\sf P}(B_{i_1,...,i_v}(n)\!\cap\! B_{\tilde i_1,...,\tilde i_v}(n))=D^2\left(\prod_{\ell=2}^{v}\frac{1}{2(i_\ell-1)}\right)\times\left(\prod_{\ell=2}^{v}\frac{1}{2(\tilde i_\ell-1)}\right)
$$
$$
\times \left(\prod_{\ell=1}^{\mu}\left( \prod_{t=\sigma_{\ell}+1}^{\sigma_{\ell+1}-1}\left(1-\frac{2\ell-1}{2(t-1)}\right)\right)\right)
\left(\prod_{\ell=\mu+1}^{2v}\left( \prod_{t=\sigma_{\ell}+1}^{\sigma_{\ell+1}-1}\left(1-\frac{2\ell-2}{2(t-1)}\right)\right)\right)
$$
$$
=\frac{D^2}{2^{2(v-1)}}\frac{\sqrt{i_1\tilde i_1}}{n^{2v-1}}\left(1+O\left(\frac{1}{i_1}\right)\right)
={\sf P}(B_{i_1,...,i_v}(n)){\sf P}(B_{\tilde i_1,...,\tilde i_v}(n))\left(1+O\left(\frac{1}{i_1}\right)\right).
$$
Therefore,
\begin{align*} 
\mathrm{Var}X&=\mathrm{Var}\left(\sum_{1\leq i_1<\ldots<i_v\leq n}I_{B_{i_1,...,i_v}(n)}\right)\\
&=\sum_{1\leq i_1<\ldots<i_v\leq n}\mathrm{Var}\left[I_{B_{i_1,...,i_v}(n)}\right]\\
&+ \sum_{{\scriptsize \begin{array}{c}1\leq i_1<\ldots<i_v\leq n,\\1\leq \tilde i_1<\ldots<\tilde i_v\leq n, \\(i_1,\ldots,i_v)\neq(\tilde i_1,\ldots,\tilde i_v)\end{array}}}
\mathrm{cov}\left(I_{B_{i_1,...,i_v}(n)},I_{B_{\tilde i_1,...,\tilde i_v}(n)}\right)\\
&<\sum_{1\leq i_1<\ldots<i_v\leq n}{\sf P}_{B_{i_1,...,i_v}(n)}\\
&+2\sum_{{\scriptsize \begin{array}{c}1\leq i_1<\ldots<i_v\leq n,\\i_1< \tilde i_1<\ldots<\tilde i_v\leq n, \\(i_1,\ldots,i_v)\cap(\tilde i_1,\ldots,\tilde i_v)=\varnothing\end{array}}} \left[{\sf P}\left(B_{i_1,...,i_v}(n)\cap B_{\tilde i_1,...,\tilde i_v}(n)\right)\right.\\
&\left.-{\sf P}\left(B_{i_1,...,i_v}(n)\right){\sf P}\left(B_{\tilde i_1,...,\tilde i_v}(n)\right)\right]\\
&<{\sf E}X+2{\sf E}X\sum_{i_1<\ldots<i_v}{\sf P}(B_{i_1,...,i_v}(n))O\left(\frac{1}{i_1}\right)
=O(n).\\
\end{align*}
Finally,
$$
 {\sf P}(X=0)\leq\frac{\mathrm{Var}X}{({\sf E}X)^2}\to 0,\quad n\to\infty.
$$


\begin{thebibliography}{99}

\bibitem{Bollobas} B. Bollob\'{a}s, {\it Random Graphs}, 2nd Edition, Cambridge University Press, 2001.

\bibitem{recursive} B. Bollob\'{a}s, O. Riordan, J. Spencer, G. Tusn\'{a}dy, {\it The degree sequence of a scale-free random graph process}. Random Structures \& Algorithms, 2001, {\bf 18}(3):~279--290. 

\bibitem{Glebsk} Y. V. Glebskii, D. I. Kogan, M. I. Liogon’kii, V. A. Talanov. {\it Range and degree of realizability of formulas in the restricted predicate calculus}. Cybernetics and Systems Analysis, 1969, {\bf 5}(2):~142--154. (Russian original: Kibernetika, 1969, {\bf 5}(2):~17--27).

\bibitem{Fagin} R. Fagin. {\it Probabilities in finite models.} J. Symbolic Logic, 1976, {\bf 41}:~50--58.

\bibitem{Haber} S. Haber, M. Krivelevich M. {\it The logic of random regular graphs}. J. Comb., 2010, {\bf 1}(3-4):~389--440.

\bibitem{Muller} P.~Heinig, T.~Muller, M.~Noy, A.~Taraz, {\it Logical limit laws for minor-closed classes of graphs}, Journal of Combinatorial Theory, Series B. 2018, {\bf 130}:~158--206.

\bibitem{Janson} S. Janson, T. Luczak, A. Rucinski, {\it Random Graphs}, New York, Wiley, 2000.

\bibitem{Kaufmann} M. Kaufmann, S. Shelah. {\it On random models of finite power and monadic logic}. Discrete Mathematics, 1985, {\bf 54}(3):~285--293.

\bibitem{Kleinberg} R. D. Kleinberg, J. M. Kleinberg. {\it Isomorphism and embedding problems for infinite limits of scale-free graphs}. In Proceedings of the 16th ACM-SIAM Symposium on Discrete Algorithms, pages 277--286, 2005.

\bibitem{Libkin} L. Libkin. {\it Elements of finite model theory}. Texts in Theoretical Computer Science. An EATCS Series. Springer-Verlag Berlin Heidelberg. 2004.

\bibitem{geo} G.L. McColm. {\it First order zero-one laws for random graphs on the circle.} Random Structures and Algorithms, {\bf 14}(3):~239--266, 1999.

\bibitem{McColm} G.L. McColm. {\it MSO zero-one laws on random labelled acyclic graphs}. Discrete Mathematics, 2002, {\bf 254}:~331--347.


\bibitem{Zhuk_Ostr_APAL} L.B. Ostrovsky, M.E. Zhukovskii. {\it Monadic second-order properties of very sparse random graphs.} Annals of pure and applied logic. 2017, {\bf 168}(11):~2087--2101.

\bibitem{Survey} A.M.~Raigorodskii, M.E.~Zhukovskii. {\it Random graphs: models and asymptotic characteristics}, Russian Mathematical Surveys, {\bf 70}(1):~33--81, 2015.

\bibitem{Spencer_Ehren} J.H. Spencer. {\it Threshold spectra via the Ehrenfeucht game.} Discrete Applied Math., 1991, {\bf 30}:~235--252.

\bibitem{Shelah} S. Shelah, J.H. Spencer. {\it Zero-one laws for sparse random graphs.} J. Amer. Math. Soc., 1988, {\bf 1}:~97--115.

\bibitem{Strange} J.H.~Spencer, {\it The Strange Logic of Random Graphs}, Springer Verlag, 2001.

\bibitem{Zhuk_Svesh} N.M. Sveshnikov,  M.E. Zhukovskii, {\it First order zero-one law for uniform random graphs}, Sbornik Mathematics, 2020, {\bf 211}, https://doi.org/10.1070/SM9321.


\bibitem{Verb} O. Verbitsky, M. Zhukovskii. {\it The Descriptive Complexity of Subgraph Isomorphism Without Numerics}, Lecture Notes in Computer Science,  International Computer Science Symposium in Russia. 2017. P.~308--322.


\bibitem{Tyszkiewicz} J.~Tyszkiewicz, {\it On Asymptotic Probabilities of Monadic Second Order Properties}, Lecture Notes in Computer Science, 1993, {\bf 702}:~425--439.

\bibitem{Shelah} S. Shelah, J.H. Spencer, \emph{Zero-one laws for sparse random graphs}, J. Amer. Math. Soc., 1988, \textbf{1}:97--115.

\bibitem{Woods} A.R. Woods, {\it Colouring rules for Knite trees, and probabilities of monadic second order sentences}, Random Structures \& Algorithms, 1997, {\bf 10}: 453--485.

\bibitem{Zhuk_JML} M.E. Zhukovskii, {\it Logical laws for short existential monadic second-order sentences about graphs}, Journal of Mathematical Logic, 2019, https://doi.org/10.1142/S0219061320500075.


\end{thebibliography}
\end{document}